\documentstyle{amsppt}
\magnification=1200
\hsize=150truemm
\vsize=224.4truemm
\hoffset=4.8truemm
\voffset=12truemm

\NoRunningHeads

\define\C{{\bold C}}
\define\R{{\bold R}}
\let\thm\proclaim
\let\fthm\endproclaim
 


\newcount\refno
\global\refno=0

\def\nextref#1{
      \global\advance\refno by 1
      \xdef#1{\the\refno}}

\def\bref {\ref\global\advance\refno by 1\key{\the\refno}}

\nextref\AH
\nextref\BR
\nextref\DE
\nextref\DT
\nextref\CU
\nextref\DG
\nextref\ST
\topmatter

\title 
Singularities of local Ahlfors currents
 \endtitle

\author  Julien Duval \footnote""{Laboratoire de Math\'ematiques d'Orsay, Univ. Paris-Sud, CNRS, Universit\'e Paris-Saclay, 91405 Orsay, France \newline julien.duval\@math.u-psud.fr\newline}
\footnote""{keywords : currents \newline AMS class. : 32C30.\newline}
\endauthor

\abstract\nofrills{\smc Abstract.} We prove that a complex curve charged by a local Ahlfors current is either a disc or an annulus.
\endabstract

\endtopmatter 

\document

\subhead 1. Introduction \endsubhead

\null
Let $B$ be the unit ball of $\C^2$. A (local) {\it Ahlfors current} is a positive current $T= \lim \frac {[D_n]}{a_n}$ where $D_n$ is a complex disc in $B$ of area $a_n$ whose boundary is mainly on $\partial B$. Precisely $\text{length}(\partial D_n\setminus \partial B)=o(a_n)$. This insures that $T$ is closed in $B$. These currents are the local version of the currents coming from an entire curve in complex hyperbolicity [\BR]. They appear naturally in geometric function theory [\DG]. We want to analyze their singular part. According to Siu (see for instance [\DE]) we may split $T$ in a sum of a diffuse part where Lelong numbers vanish except at a countable number of points and a singular part consisting in a combination of currents of integration on complex curves in $B$. These curves are precisely those charged by $T$.  
Here is our result.
 
\thm{Theorem} Let $T$ be an Ahlfors current in $B$ and $C$ a complex irreducible curve charged by it. Then $C$ is either a disc or an annulus.
\fthm
 
This result is the local counterpart of a similar statement which says that an irreducible algebraic curve in $P^2(\C)$ charged by a (global) Ahlfors current is either rational or elliptic [\CU]. Both proofs share the same scheme but differ in a crucial point, the way to achieve the lifting property of paths from $C$ to $D_n$. In the global case this goes by using Ahlfors theory of covering surfaces which cannot readily apply in the local situation because of the presence of a large boundary in $\partial B$. Here we are able to lift very short paths thanks to the laminar structure of $T$ [\DT]. The new point is to analyze the cost in terms of area of the failure to extend this lifting along larger paths. 

Let us sketch the proof of our theorem. We argue by contradiction. Then there exists an essential figure eight $\gamma$ in $C$. According to this lifting property we may lift $\gamma$ in $D_n$ in a nearby graph $\gamma_n$ with many cycles. This reflects the fact that the projection of a tubular neighborhood to $C$ when restricted to $D_n$ behaves very much like a covering map of large degree. Hence we get a lot of discs in $D_n \setminus \gamma_n$. Now on one hand such a disc is close to $\gamma$ by the maximum principle. On the other hand its boundary is close to a non trivial loop in $\gamma$. This is the contradiction.

Let us detail this.

\null

\subhead 2. Proof of the theorem \endsubhead

\null
Let $T= \lim \frac {[D_n]}{a_n}$ be an Ahlfors current in the unit ball $B$ of $\C^2$. Again $D_n$ is a complex disc in $B$, smooth up to the boundary, of area $a_n$ and such that $\text{length}(\partial D_n\setminus \partial B)=o(a_n)$. We suppose that $T$ charges an irreducible complex curve $C$ of $B$. So there exists $\nu >0$ such that $T\vert_C=\nu[C]$. Up to shrinking a bit $B$ we may assume that $C$ is a compact Riemann surface with smooth boundary.

\null\noindent
{\bf The argument}. 

\null
 We argue by contradiction supposing that $C$ is neither a disc nor an annulus. Then we may find an {\it essential} figure eight $\gamma$ in its smooth part. Here essential means that all components of $C\setminus \gamma$ reach $\partial C$. It follows that $\gamma$ is polynomially convex. We refer to [\ST] for background on polynomial convexity. Call $E$ a thin (still polynomially convex) neighborhood of $\gamma$ in $C$. Trivializing the normal bundle of $C$ on $E$ we get an ambiant polynomially convex neighborhood $V$ of $\gamma$ of the form $E\times D$ where $D$ is a small disc. Call $\pi$ the first projection. Note that $d_n=V\cap D_n$ consists in a finite union of complex discs because of the polynomial convexity of $V$. We will lift (a small perturbation of) $\gamma$ through $\pi_n=\pi\vert_{d_n}$
in a graph $\gamma_n$ with a lot of cycles. This is enough to reach a contradiction. Indeed by a simple topological argument such a graph has to cut out a lot of discs in $d_n$. Take one of them. By construction it projects down by $\pi$ to a complex disc in $E$ with (non constant) boundary in $\gamma$. But by polynomial convexity of $\gamma$ this projection has to be contained in $\gamma$. Therefore it is constant.

Let us turn now to the main point.

\null\noindent
{\bf Lifting paths}.

\null
Fix a small $\epsilon>0$ (say $<\frac \nu 9$). First we may find many lifts near a generic point of $E$ because of the laminarity of $T$. We refer to [\DT] for background on laminar currents and details on this fact. In particular we may fix a small disc $d$ in $E$ containing (after perturbation) the vertex of $\gamma$ such that $\pi_n^{-1}(d)$ contains $(\nu-\epsilon)a_n$ discs which are graphs over $d$ and converge to $d$ when $n \to \infty$. We want to extend these lifts along the two loops of $\gamma$. Consider one of them given by an arc $\alpha$ starting and ending in $d$. Fix a large integer $k$ (say $> \frac {2\nu} \epsilon$). Take $k$ disjoint thin strips parallel to $\alpha$ in $E$, each of them foliated by arcs parallel to $\alpha$. In the sequel we will ignore a countable number of these arcs, the {\it critical} ones (those containing a critical value of $\pi_n$). Near the starting points of the arcs we already have $(\nu-\epsilon)a_n$ lifts. 

We will see now that we can extend most of them along the arcs of one strip. Actually we will estimate the cost of the failure of this extension. For this consider a lift of $d$ and suppose that for 3 strips there exists an arc along which the initial lift cannot be fully extended. Note that we still get partial lifts. Call such a lift of $d$ together with these 3 partial lifts a {\it tripod}. It sits in $d_n$ with 3 end points in $\partial d_n$. The point is the following

\thm{Fact} The number of tripods is negligible with respect to $a_n$.
\fthm

Suppose not. Then we have plenty of (disjoint) tripods, say $\eta a_n$ (for some $\eta>0$). Their union cuts out plenty of discs in $d_n$. We are interested in the {\it terminal} discs, those which touch a single tripod. A simple combinatorial argument tells us that they are in number $\geq \eta a_n$. Indeed consider a t-vertex for each tripod, a d-vertex for each disc in the complement of the tripods and an edge between them if such a disc touches the tripod. In this way we construct trees. Call $v_t, v_d$ the number of t- or d-vertices and $i$ the number of incidence edges. We have $i=3v_t$ (because a tripod touches 3 discs) and $v_t+v_d\geq i$ (because we have trees) so that $v_d\geq 2v_t$. Call $v_{td}$ the number of terminal discs. Counting the edges in a different way we get $i\geq v_{td}+2(v_d-v_{td})$. Therefore $v_{td}\geq 2v_d-3v_t\geq v_t$ which concludes.

But this contradicts the following lemma (whose proof we postpone). 
\thm{Lemma} Let $\delta_n$ be a sequence of terminal discs such that length($\partial \delta_n\cap \partial D_n )\to 0$.
Then area($\delta_n) \to \infty$.
\fthm
Indeed order the terminal discs $\delta$ by increasing length of $\partial \delta\cap \partial D_n$ and retain only the first half of them. We then have $$ \text{length}(\partial \delta\cap \partial D_n) \leq \frac {\text{2 length} (\partial D_n \setminus \partial B)}{\eta a_n}$$ which goes to $0$ by assumption. So we end up with $\frac \eta 2 a_n$ disjoint discs in $D_n$ of area going uniformly to infinity. This is impossible as the total area of $D_n$ is $a_n$. 

\null

At this stage we know that for $(\nu-2\epsilon) a_n$ of them, we are able to extend a lift of $d$ along all the arcs of $k-2$ strips. By a simple pigeonhole argument we choose a good strip $\sigma$ such that $(\nu-3\epsilon)a_n$ lifts of $d$ extend along all the arcs of $\sigma$.
Shrinking $V$ toward $E$ and repeating the process, we may arrange by a diagonal argument that these lifts converge uniformly toward $\sigma$.

\newpage\noindent
{\bf Gluing lifts}.

\null
We now check that most of these lifts along a good arc of $\sigma$ glue back at the end to the lifts of $d$.
Indeed denote by $q$ the component of $\sigma \cap d$ corresponding to the end part of $\sigma$. It is a small square of size $r$ foliated in one direction by arcs $\alpha_t$ parallel to the end part of $\alpha$. For each $\alpha_t$  we know of 2 types of lifts : the first are the restriction to $\alpha_t$ of the lifts of $d$, the second the continuation through $\sigma$ of the lifts of $d$. These various lifts either are disjoint or coincide. Denote by $s(t)$ the number of lifts of $\alpha_t$ of the second type but not of the first. We estimate $s(t)$ by an area argument. Note on one hand that the area of the part of $d_n$ close to $q$ is $\leq (\nu+\epsilon)a_n r^2$ as $T\vert_C=\nu[C]$. On the other hand lifts of the first type already fill out an area $\geq(\nu-\epsilon) a_n r^2$. Estimating the area for the rest of the lifts we get

 $$\int_0^rrs(t)dt \leq 2 \epsilon a_n r^2 \text{ so }\min_{t\in [0,r]} s(t)\leq 2 \epsilon a_n.$$
 This singles out our good arc.
Dealing in the same way with the other loop of $\gamma$ we end up with a graph $\gamma_n$ in $d_n$ with $(\nu-\epsilon)a_n$ vertices (the lifts of $d$) and $2(\nu - 5\epsilon)a_n$ edges (the lifts of the 2 good arcs, starting and ending at a lift of $d$). Its Euler characteristic satisfies $$\chi(\gamma_n)\leq -(\nu-9\epsilon)a_n.$$
We do have plenty of cycles in $\gamma_n$ as expected.

\null\noindent
{\bf Proof of the lemma}.

\null

We keep the notations above.
We argue by contradiction supposing the area of the terminal discs $\delta_n$ bounded by some constant $C$. Note that by construction $\partial \delta_n$ splits in 4 parts, 3 of them being in the tripod (a part of a lift of $d$ and 2 partial lifts of arcs) and the last in $\partial d_n$. Parametrize $\delta_n$ by a rectangle accordingly. We get $f_n : R_n=[0,m_n]\times [0,1](\subset \R^2\simeq \C) \to \delta_n$, holomorphic and smooth up to the boundary, sending
the horizontal parts of $\partial R_n$ to the partial lifts, the left part in the lift of $d$ and the right in $\partial d_n$. This determines the modulus $m_n$. We refer to [\AH] for background in conformal geometry. 

We first note that $m_n$ remains bounded. Indeed if not we would find by a lenght-area argument [\AH] a vertical segment $I_n$ in $R_n$ such that $\text{length}(f_n(I_n))\to 0$. But this is impossible as $f_n(I_n)$ connects the 2 partial lifts which remain far apart because they project down in 2 different strips. After extraction we may think the modulus (almost) constant.

In a same vein $f_n(\{m\}\times[0,1])$ is close to $\partial V$. Indeed it is included in $\partial d_n=(\partial D_n \cap V)\cup (D_n \cap \partial V)$. By assumption we have length($\partial \delta_n\cap \partial D_n )\to 0$. Again as the 2 partial lifts are far apart this forces $f_n(\{m\}\times[0,1])$ to touch $\partial V$, and to remain close to it.

We also know that $f_n(\{0\}\times [0,1])$ converges toward $E$ (it is part of a lift of $d$). This implies that $f_n$ converges locally uniformly toward $E$ in the interior of $R$ (apply the two-constant theorem [\AH] to $ p\circ f_n $ where $p$ is the second projection of $V=E\times D$).

The contradiction will be reached at the right bottom corner $b$ of $R$ by a length-area argument which we detail. Take polar coordinates centered at $b$ and denote by $\lambda_n \vert dz \vert$ the pull-back of the euclidean metric by $f_n$. We are looking for a quarter of circle centered at $b$ of controlled radius but short image. Its length is given by $l_n(r)=\int_{\frac \pi 2}^\pi\lambda_n rd\theta$. Estimating the area of the image of a quarter of annulus of radii $\rho\ll \rho'$ we get 
$$\int_{\rho}^{\rho'}\int_{\frac \pi 2}^\pi{\lambda_n}^2 rdrd\theta\leq C \hskip 3mm \text{so }\min_{r\in [\rho,\rho']} l_n(r)\leq \sqrt{\frac{\pi C}{2\log \frac{ \rho'}{\rho}}}$$
by Cauchy-Schwarz inequality.  
This singles out our good quarter of circle. Denote by $x_n,y_n,z_n$ the images of its middle, its lower end and its right end. They are close to each other. By the discussion above $x_n$ is close to $E$ and $z_n$ to $\partial V$. By construction $y_n$ projects down in a given strip.  This forces this strip to be (arbitrarily) close to $\partial V$, which is impossible.

\Refs
\widestnumber\no{99} 
\refno=0 
\bref \by L. Ahlfors \book Conformal invariants. Topics in geometric function theory  \publ AMS Chelsea Publ. \yr2010 \publaddr Providence
\endref
 
\bref \by M. Brunella \paper Courbes enti\`eres et feuilletages holomorphes \jour Enseign. Math. \vol45\yr1999\pages195--216
\endref
\bref \by J.-P. Demailly \book Monge-Amp\`ere operators, Lelong numbers and intersection theory,  {\rm in}
 Complex analysis and geometry \pages115--193 \publ Plenum \yr1993 \publaddr New-York
\endref
\bref \by H. de Th\'elin \paper Sur la laminarit\'e de certains courants \jour Ann. Sci. ENS \vol37\yr2004\pages304--311
\endref
\bref \by J. Duval \paper Singularit\'es des courants d'Ahlfors \jour Ann. Sci. ENS \vol39\yr2006\pages527--533
\endref
\bref \by J. Duval and D. Gayet \paper Riemann surfaces and totally real tori \jour Comment. Math. Helv. \vol89\yr2014\pages299--312
\endref
\bref \by E. L. Stout \book Polynomial convexity \bookinfo Prog. Math. \vol261 \publ Birkh\"auser \yr2007 \publaddr Boston
\endref
 
\endRefs

\enddocument